\newcommand{\frr}[1]{}

\def\Z{{\mathbb Z}}
\def\Q{{\mathbb Q}}
\def\R{{\mathbb R}}

\def\Zt{\Z[\tau]}
\def\Zb{\Z[\beta]}

\def\Qb{\Q[\beta]}

\def\pf{\begin{proof}}
\def\pfk{\end{proof}}

\documentclass{amsart}
\usepackage{epsfig}
\begin{document}

\newsymbol\Vdash 130D
\newsymbol\supsetneqq 2325
\newsymbol\subsetneqq 2324
\newtheorem{lem}{Lemma}[section]
\newtheorem{thm}[lem]{Theorem}
\newtheorem{prop}[lem]{Proposition}
\newtheorem{coro}[lem]{Corollary}
\newtheorem{que}{Question}

\theoremstyle{definition}
\newtheorem{de}[lem]{\bf Definition}

\theoremstyle{remark}
\newtheorem{pozn}[lem]{Remark}


\title{$s$-convexity, model sets and their relation}

\author{Zuzana Mas\'akov\'a}
\address{Department of Mathematics,
Faculty of Nuclear Science and Physical Engineering\\
Czech Technical University, Trojanova 13, Prague 2, 120 00, Czech Republic\\
masakova@km1.fjfi.cvut.cz}

\author{Ji\v r\'\i\  Patera}
\address{Centre de recherches math\'ematiques,
Universit\'e de Montr\'eal, \\
C.P. 6128 Succursale Centre-ville, Montr\'eal, Qu\'ebec,
H3C 3J7, Canada\\ 
patera@crm.umontreal.ca}

\author{Edita Pelantov\'a}
\address{Department of Mathematics,
Faculty of Nuclear Science and Physical Engineering\\
Czech Technical University, Trojanova 13, Prague 2, 120 00, Czech Republic\\
pelantova@km1.fjfi.cvut.cz}


\begin{abstract}
The relation of $s$-convexity and sets modeling physical quasicrystals 
is explained for quasicrystals related to quadratic unitary Pisot numbers. 
We show that 1-dimensional model sets may be characterized by
$s$-convexity for finite set of parameters $s$. It is shown that the
three Pisot numbers $\frac12(1+\sqrt5)$, $1+\sqrt2$, and $2+\sqrt3$
related to experimentally observed non-crystallographic symmetries
are exceptional with respect to $s$-convexity.
\end{abstract}

\maketitle

\section{Introduction}
\label{sec:prelim}

In this paper we study properties of a class of binary operations defined
on a real vector space.
As we show in this article such operations can be used
for characterizations of point sets modeling physical quasicrystals. 
Such operation was first introduced in a purely mathematical article 
by I.~Calvert in 1978.
However, from his paper~\cite{calvert} and absence of references in it, 
one cannot deduce the motivation for the study of this operations. The
problem was then elaborated by R.~G.~E.~Pinch in 1985 in~\cite{pinch}. 
The connection of this operation with quasicrystals was first recognized 
by Berman and Moody in~\cite{bm}, where they work with a special case 
of such operation. 

For arbitrary real parameter $s$, Calvert defined a binary operation
\begin{equation}\label{eq:oper}
x\vdash_s y:=sx+(1-s)y\,, \qquad x,y\in\R^n\,.
\end{equation}
Pinch in~\cite{pinch} calls a set closed under such operation an 
{\bf $s$-convex} set. Pinch shows that an $s$-convex set $\Lambda\subset\R$ 
is either dense in an interval of $\R$,
or it is uniformly discrete. A set $\Lambda\subset\R$ is called 
{\bf uniformly discrete}, if there exists an $\varepsilon >0$, 
such that $|x-y|>\varepsilon$ for any $x,y\in\Lambda$, $x\neq y$. 
Pinch studies the question for which parameters $s$ there exists 
a uniformly discrete $s$-convex set containing at least two points.
We shall denote by ${\it Cl}_sA$ the closure of the set $A$ under the 
operation~\eqref{eq:oper}.
A most suitable candidate for a uniformly discrete $s$-convex set 
($s$ fixed) is ${\it Cl}_s\{0,1\}$.
If $s$ belongs to $(0,1)$, then ${\it Cl}_s\{0,1\}$ is a closure under 
a ordinary convex combination and therefore is dense in $[0,1]$. 
Pinch shows that for $s\notin(0,1)$, the closure ${\it Cl}_s\{0,1\}$ 
may be uniformly discrete only if $s$ is an algebraic integer. On 
the other hand he proves that if $s$ is a totally real algebraic integer, 
such that all its algebraic conjugates are in $(0,1)$, then ${\it Cl}_s\{0,1\}$ 
is uniformly discrete. 

The connection to sets modeling quasicrystals was observed by 
Berman and Moody on the example of model sets based on the 
golden ratio $\tau=\frac12(1+\sqrt5)$, defined by
\begin{equation}\label{eq:tau}
\Sigma_\tau(\Omega) := \{a+b\tau \mid a,b\in\Z\,,\, a+b\tau' \in \Omega\}\,,
\end{equation}
where $\tau'=\frac12(1-\sqrt5)$ is the algebraic conjugate of $\tau$ and
$\Omega$ is a bounded interval in $\R$. For such $\Omega$, 
the set  $\Sigma_\tau(\Omega)$
is not only uniformly discrete, but also {\bf relatively dense}, which means 
that the distances of adjacent points in $\Sigma_\tau(\Omega)$ are
bounded by a fixed constant. A set which is both uniformly discrete
and relative dense is called a {\bf Delone} set.
It is easy to show that $\Sigma_\tau(\Omega)$ is closed under the 
operation~\eqref{eq:oper} for $s=-\tau$. 
The operation $\vdash_{-\tau}$ has an important 
property. It has been shown in~\cite{bm} that the quasicrystal
$\Sigma_\tau[0,1]$ may be generated by this operation starting 
from initial seeds $\{0,1\}$, more precisely ${\it Cl}_{-\tau}\{0,1\}=
\Sigma_\tau[0,1]$. This means that the operation is `strong enough'
so that any uniformly discrete set closed under it has to be a model set
of the form~\eqref{eq:tau}.

It can be shown that all model sets~\eqref{eq:tau} are $s$-convex 
for infinitely many $s\in\Zt$.  
A question arises whether there exists another $s$ `strong enough'
to ensure ${\it Cl}_{s}\{0,1\}=\Sigma_\tau[0,1]$. In this article we provide
an answer to a more general question. 

The golden ratio $\tau$ is a quadratic unitary Pisot number related 
to 5-fold symmetry quasicrystals. Other experimentally observed 
non-crystallographic symmetries (8-fold and 12-fold) are related
also to quadratic unitary Pisot numbers $1+\sqrt2$ and $2+\sqrt3$.
We consider here two infinite families of quadratic unitary Pisot numbers
$\beta>1$, solutions of the equations
\begin{eqnarray*}
&&x^2=mx+1\,, \qquad m\in\Z\,, m\geq1\,, \hbox{ or}\\
&&x^2=mx-1\,, \qquad m\in\Z\,, m\geq3\,.
\end{eqnarray*}
Let us denote by $\beta'$ the algebraic conjugate of $\beta$, i.e.~the 
second root of the equation. The number $\beta$ is called unitary 
Pisot, since $|\beta'|<1$, and $\beta\beta'=\mp1$.
A model set is then defined as a subset of the ring $\Z[\beta]:=\Z+\Z\beta$,
equipped with the Galois automorphism $x=a+b\beta$ $\ \mapsto\ $ 
$x'=a+b\beta'$. We let
\begin{equation}\label{eq:deqc}
\Sigma_\beta(\Omega) := \{x\in\Zb \mid x'\in\Omega\}\,,
\end{equation}
where $\Omega$ is a bounded interval in $\R$. 
The set $\Sigma_\beta(\Omega)$ is called a model set in $\Zb$ with 
the acceptance window $\Omega$. 

Remark: Let us mention that the solution of equation $x^2=3x-1$ gives
the same ring $\Zb$ as the equation $x^2=x+1$, and thus we shall consider
the equation $x^2=mx-1$ with $m>3$ only.  

\begin{que}\label{q:1}
Let $\beta$ be quadratic unitary Pisot number. Which parameter 
$s$ ensures that any uniformly discrete $s$-convex set
is a model set in $\Zb$? 
\end{que}

We show in this paper that $s=-\frac12(1+\sqrt5)$ is exceptional 
for model sets in $\Zt$. The number $\tau^2=1-s$ has the same 
property, since generally $s$-convexity is equivallent to 
$(1-s)$-convexity. For no other parameters $s$ the 
operation~\eqref{eq:oper} is `strong enough'. 
The only other quadratic unitary Pisot numbers, for which a 
suitable parameter exists are $\beta=1+\sqrt2$ and $2+\sqrt3$. 
It means that in the majority of cases the invariance with respect 
to one operation~\eqref{eq:oper} does not imply being a model set. 
One may further ask the following question.

\begin{que}\label{q:2}
Let $\beta$ be quadratic unitary Pisot number. 
Is it possible to find a finite set ${\mathcal N}$ of parameters,
such that $s$-convexity of a uniformly discrete set $\Lambda\subset\Zb$
for all $s\in {\mathcal N}$ implies $\Lambda$ being a model set in $\Zb$?
\end{que}

The answer is YES, and we provide a constructive proof for each $\beta$.


The solutions to Questions~\ref{q:1} and~\ref{q:2} are formulated as
Theorems~\ref{thm:1} and~\ref{thm:2} in Section~\ref{sec:results}.
Their proofs require different approaches and different mathematical 
tools, therefore they are found in separate Sections~\ref{sec:tools} 
and~\ref{sec:druha}.

\section{Auxiliary lemmas I}\label{sec:tools}

Recall that we are working in a ring $\Zb\subset\Qb$, and that the
mentioned conjugation $\ '\ $ is a Galois automorphism on the field
$\Qb$. One has $(x+y)'=x'+y'$ and $(xy)'=x'y'$ for any two elements
$x,y\in\Zb$. With this in mind, we may state that for $s\in\Zb$ and  
$\Lambda\subset\Zb$,  
\begin{equation}\label{eq:eqv}
\Lambda \hbox{ is } s\hbox{-convex} \qquad\hbox{ if and only if }\qquad
\Lambda'=\{x'\mid x\in\Lambda\} \hbox{ is } s'\hbox{-convex.}
\end{equation}
From the definition of model sets it follows that $\bigl(\Sigma_\beta
(\Omega)\bigr)'=\Zb\cap\Omega$. Hence from~\eqref{eq:eqv} we obtain
$$
\Sigma_\beta(\Omega) \hbox{ is } s\hbox{-convex} \qquad\hbox{ if and only if }
\qquad \Zb\cap\Omega \hbox{ is } s'\hbox{-convex.}
$$
Since $\Omega$ is an interval, $s'x+(1-s')y$ has to be a convex 
combination in the ordinary sense, i.e.~$s'\in[0,1]$, which is equivallent to
$s\in\Sigma_\beta[0,1]$.

In this formalism, Question~\ref{q:1} may be rewritten as follows:
Does there exist $s\in\Sigma_\beta[0,1]$, such that for any uniformly 
discrete $s$-convex set $\Lambda\subset\Zb$ it holds that
\begin{equation}\label{eq:ladnenko}
\Lambda'=\Zb\cap\Omega
\end{equation}
for some $\Omega$?
If~\eqref{eq:ladnenko} should be valid for every uniformly discrete $\Lambda$, 
it must be true also for $\Lambda={\it Cl}_{s}\{0,1\}$. Since $\Lambda'=
{\it Cl}_{s'}\{0,1\}$ contains $0'=0$ and $1'=1$ and since $s'x+(1-s')y$ is
a convex combination, the only suitable candidate for $\Omega$ in this 
case is the interval $\Omega=[0,1]$. It is the reason why in the sequel
we focus our attention to investigation of equatlity
$$
[0,1]\cap\Zb={\it Cl}_{s'}\{0,1\}\,, \qquad \hbox{ for } s'\in(0,1)\,.
$$

\begin{lem}\label{xorxj}
Let $\beta$ be quadratic unitary Pisot number. Assume that
${\it Cl}_{s'}\{0,1\}=\Zb\cap[0,1]$. Then for any $y\in\Zb\cap[0,1]$, the
scaling factor $s'$ divides either $y$ or $y-1$ in the ring $\Zb$.
\end{lem}

\pf
Pinch in~\cite{pinch} has proved that any $y\in{\it Cl}_{s'}\{0,1\}$ may be
written in the form 
\begin{equation}\label{eq:as}
y=\sum_{i=0}^n b_i(s')^i(1-s')^{n-i}\,, \qquad b_i\in\Z\,,\, 
0\leq b_i \leq {n\choose i}\,,
\end{equation}
for some non negative integer $n$. 
Therefore $b_0\in\{0,1\}$. 
If $b_0$ is equal to 0, then $y$ is divisible by $s'$, otherwise $s'$ divides
$y-b_0=y-1$.
\pfk

Before the following corollary let us mention several number theoretical facts.
On the field $\Qb$ on may define a `norm' $N(x):=xx'\in\Q$ and a `trace'
${\rm tr}(x):=x+x'\in\Q$. Since $|N(\beta)|=1$ and ${\rm tr}(\beta)=m$ are both
integers, it holds for every $x\in\Zb$ that $N(x)\in\Z$ and ${\rm tr}(x)\in\Z$.
A divisor of unity in $\Zb$ is an element $u$ such that $\frac1u\in\Zb$; $u$
is a divisor of unity iff $|N(u)|=1$.

\begin{coro}\label{divdva}
If ${\it Cl}_{s'}\{0,1\}=\Zb\cap[0,1]$, then both $s'$ and $1-s'$
divide 2 in the ring $\Zb$.
\end{coro}

\pf
Note that the roles of $s'$ and $1-s'$ are symmetric. Thus it suffices to 
show that $s'$ divides $2$. 
If $s'$ is a divisor of unity then the assertion is true.
Assume the opposite, i.e. $s'$ is not a divisor 
of unity and $s'$ does not divide $2$. Choose $y=1/\beta^2$. 
Since $y$ is a divisor of unity, 
the scaling factor $s'$ does not divide neither $y$, 
nor $2y$. According to Lemma~\ref{xorxj} this
implies $s'|(y-1)$ and $s'|(2y-1)$. Therefore $s'$ divides
$(2y-1)-(y-1)=y$, thus a contradiction.    
\pfk


\begin{prop}\label{vylouc}
Let $\beta$ be a quadratic unitary Pisot number. Assume that
${\it Cl}_{s'}\{0,1\}=\Zb\cap[0,1]$ is satisfied for some $s'$. 
Then either of the possibilities below is true:
\begin{itemize}
\item $\beta=\tau$ (root of $x^2=x+1$) and $s'=-\tau'$ or $s'=1+\tau'$.
\item $\beta=1+\sqrt2$ (root of $x^2=2x+1$) and $s'=-\beta'$ or $s'=1+\beta'$.
\item $\beta=2+\sqrt3$ (root of $x^2=4x-1$) and $s'=\beta'$ or $s'=1-\beta'$.
\end{itemize}
\end{prop}

\pf
From the relation~\eqref{eq:as} it is clear that all elements of ${\it Cl}_{s'}\{0,1\}$
are polynomials in $s'$. Since $s'\in\Zb$, it satisfies a quadratic equation
with integer coefficients, namely $x^2=(s+s')x-ss'$, i.e. $x^2={\rm tr(s)}x-N(s)$. 
Using this quadratic equation any polynomial $y$ from~\eqref{eq:as} 
may be reduced to the form $y=a+bs'$. Therefore clearly ${\it Cl}_{s'}\{0,1\}
\subset\Z[s']$. The condition ${\it Cl}_{s'}\{0,1\}=\Zb\cap[0,1]$ implies
that we need $\Zb=\Z[s']$. This can be satisfied only if $s'=\pm\beta+k$ for
some integer $k$. The restriction $0<s'<1$ gives only two admissible values 
for the parameter $s'$. 

(a) Let us consider $\beta^2=m\beta+1$, i.e. $[\beta]=m$. Here
$s'=\beta-m=1/\beta$ or $s'=m+1-\beta=1-1/\beta$. 
Corollary~\ref{divdva} states that both $s'$ and $1-s'$ divide 2 in $\Zb$.
Since $N(1/\beta)=-1$, $s'=1/\beta$ is a divisor of unity and hence divides 
2 automatically. Let us assume that $1-1/\beta$ divides 2. It means that there
exists an element $c+d\beta\in\Zb$, such that
$$
2=(c+d\beta)\left(1-\frac1\beta\right)=
c-d +mc+\beta(d-c)\,.
$$
It follows that $c=d$ and $mc=2$. Since $m$ is a positive integer, this may
happen only for $m=1$ and $m=2$. These are the two cases given in the
statement.

(b) Let us now consider $\beta^2=m\beta-1$, i.e. $[\beta]=m-1$. Here
$s'=m-\beta=1/\beta$ or $s'=\beta-m+1=1-1/\beta$. Again from 
Corollary~\ref{divdva}, $1-1/\beta$ divides 2. It means that there
exists an element $c+d\beta\in\Zb$, such that
$$
2=(c+d\beta)\left(1-\frac1\beta\right)=c-d -mc+\beta(d+c)\,.
$$
This implies $c=-d$ and $(2-m)c=2$. For the equation $\beta^2=m\beta-1$
we consider $m\geq4$, so that the only solution is $m=4$. 
\pfk

\section{Auxiliary lemmas II}\label{sec:druha}

In order to prove the results concerning the second question, we shall
need the notion of $\beta$-expansions.

The ring $\Zb:=\Z+\Z\beta$, can be characterized 
using $\beta$-expansions~\cite{renyi}.
A $\beta$-expansion of a real number $x\geq0$ is 
defined for any real $\beta>1$ as an infinite sequence
$(x_i)_{k\geq i>-\infty}$ given by the `greedy' algorithm in the following way:
$$
x_k:=\left[\frac{x}{\beta^k}\right]\,,
$$
where $k$ satisfies $\beta^k \leq x < \beta^{k+1}$. Denote $r_k=x/\beta^k 
- x_k$. Numbers $x_{i-1}$ and $r_{i-1}$ are computed from $x_i$ and $r_i$ by
prescription:
$$
x_{i-1}:= [\beta r_i]\,, \qquad r_{i-1} := \beta r_i - x_{i-1}\,.
$$
Clearly $x_i\in\{0,1,\dots,[\beta]\}$ for each $i$. In~\cite{renyi} it is proven 
that 
$$
x=\sum_{i=-\infty}^k x_i\beta^i\,,
$$
i.e.~the sum converges for each positive real $x$ and for each $\beta>1$.

Parry in~\cite{parry} answered the question for which sequences 
$(x_i)_{k\geq i>-\infty}$ there exists a positive real $x$, such that 
$(x_i)_{k\geq i>-\infty}$ is its $\beta$-expansion. Let the R\'enyi 
$\beta$-representation of 1 be
$$
1=\frac{a_1}\beta + \frac{a_2}{\beta^2} + \frac{a_3}{\beta^3} + \cdots\,,
\qquad a_i\in\{0,1,\dots,[\beta]\}\,,
$$
and let
$$
x=\sum_{i=-\infty}^k x_i\beta^i\,, \qquad x_i\in\{0,1,\dots,[\beta]\}\,.
$$
Then the sequence $(x_i)_{k\geq i>-\infty}$ is a $\beta$-expansion of $x$ 
if and only if  for any integer $j\leq k$, the sequence $x_jx_{j-1}x_{j-2}\dots$ 
is lexicographically strictly smaller than sequence $a_1a_2a_3\dots$.
Let us apply the rule on $\beta$-expansions for quadratic unitary Pisot
numbers $\beta$.

At first let us discuss the case $\beta^2=m\beta+1$, i.e.~$[\beta]=m$.
The R\'enyi representation of 1 is $1=m/\beta + 1/\beta^2$, so that 
$a_1=m$, $a_2=1$, $a_3=0$, $\dots$. It means that $x=\sum x_i\beta^i$ is a
$\beta$-expansion iff any $x_i=m$ occuring in the sequence 
$(x_i)_{k\geq i>-\infty}$ is followed by $x_{i-1}=0$.

Let now $\beta^2=m\beta-1$, so that $[\beta]=m-1$. The R\'enyi representation
of 1 in this case is
$$
1=(m-1)/\beta + (m-2)\sum_{k=2}^{\infty} \frac1{\beta^k}\,.
$$
Therefore  $(x_i)_{k\geq i>-\infty}$ is a $\beta$-expansion of some $x>0$ iff
$x_ix_{i-1}x_{i-2}\dots$ is strictly lexicographically smaller than
$(m-1)(m-2)(m-2)\dots$ for any $i\leq k$.

Let $x\in\R$. If $x$ is negative we put $x_i =-|x|_i$, negatives of
$\beta$-expansion coefficients for $|x|$.
If the $\beta$-expansion of an $x\in\R$ ends in infinitely many zeros, 
it is said to be finite, and the zeros at the end are omitted. 
Denote the set 
$$
{\rm Fin}(\beta) := \{ \varepsilon x \mid x\in\R^+_0\,,\, \varepsilon=\pm1\,,\,
x \text{ has a finite $\beta$-expansion}\} \,.
$$
In the sequel, it will be useful to use the relation between the set 
${\rm Fin}(\beta)$ of all numbers with finite $\beta$-expansion and 
the ring $\Zb$. In~\cite{froug} it is proved that for $\beta$, which is 
a root of the equation $x^2=mx+1$, the set ${\rm Fin}(\beta)$ is a ring. 
It follows immediately that 
\begin{equation}\label{eq:burdia2}
\Zb = {\rm Fin}(\beta)\,.
\end{equation}
For $\beta$ being a solution of $x^2=mx-1$, the situation is different. 
In this case, ${\rm Fin}(\beta)$ is not closed under addition, (for example 
$\beta-1\not\in {\rm Fin}(\beta)$), while $\Zb$ is. Burd\'\i k et al. 
in~\cite{burda} prove that for such $\beta$, the set ${\rm Fin}(\beta)$ 
can be characterized as the set of those elements $x$ of $\Zb$, for 
which $N(x)=xx'$ is non negative.
Therefore
\begin{equation}\label{eq:burdia}
\Zb \supsetneqq {\rm Fin}(\beta)=\{x\in\Zb \mid xx'\geq0\}\,.
\end{equation}

The crucial proposition for the answer to Question~\ref{q:2} is the following 
one. 

\begin{prop} \label{jajaj}
Let $\beta$ be a quadratic unitary Pisot number, satisfying the equation
$x^2=mx\pm1$. Let $\Lambda$ be a uniformly discrete subset of $\Zb$, 
containing $0,1$. Put
$$
{\mathcal M}=\left\{ \frac{i}\beta \biggm| i=1,2,\dots,\left[\frac{m\pm1}2\right] 
\right\}\,.
$$
If $\Lambda'$ is $s$-convex for every $s\in {\mathcal M}$, then there exists
a bounded interval $\Omega$, such that 
$\Lambda'=\Omega\cap\Zb$.
\end{prop}  

The proof of the above proposition is based on the following lemma. 
It is useful to introduce the notion of a closure under a set of 
operations $\vdash_s$. We denote by ${\it Cl}_{\mathcal M}A$ the smallest 
set containing $A$, which is $s$-convex for all $s$ in a finite set 
${\mathcal M}$. 

\begin{lem} \label{eins}
Let $\beta$ be a quadratic unitary Pisot number and ${\mathcal M}$ defined 
in Proposition~\ref{jajaj}. Then
$$
{\it Cl}_{\mathcal M}\{0,1\}=[0,1]\cap\Zb
$$
\end{lem}

\pf
Consider the operations $\vdash_s$ with $s=i/\beta$. For simplicity, denote
$$
x\Vdash_i y := \frac{i}\beta x + \left(1-\frac{i}\beta\right)y\,,\qquad i=1,2,\dots
$$
The statement of the lemma is equivalent to the fact that any point $\Zb$ in the
interval $[0,1]$ can be generated from $0$, and $1$ using operations
$\Vdash_i $, $i=1,2,\dots,[(m\pm1)/2]$. We first show that
$[0,1]\cap\Zb$ can be generated using all operations $\Vdash_1,\dots,\Vdash_r$,
with $r=[\beta]$. In the second step we find expression for 
$\Vdash_i$, $i=[(m\pm1)/2]+1,\dots,[\beta]$,
as a composition of operations $\Vdash_i$, $i=1,\dots,[(m\pm1)/2]$.

\begin{trivlist}
\item{(a)}
Consider $\beta$ to be the solution of the equation $x^2=mx+1$, $m\geq1$.
According to~\eqref{eq:burdia2}, any $x\in\Zb\cap(0,1)$ has a finite 
$\beta$-expansion 
$$
x=\sum_{i=1}^k\frac{x_i}{\beta^i}\,, \qquad k\in\Z\,.
$$
By induction on $k$ we show that $x$ belongs to ${\it Cl}_{\mathcal M}\{0,1\}$. 
Clearly, for the first step of induction we have 
\begin{equation} \label{eq:1step}
\frac{j}\beta=1\Vdash_j 0\,,
\end{equation}
for any $j=1,\dots,m$. Suppose that the $\beta$-expansion of a point
$x\in\Zb\cap(0,1)$ is of the form
$$
x=\frac{j}\beta + \sum_{i=2}^k\frac{x_i}{\beta^i}\,, \qquad j=0,\dots,m-1\,.
$$
Then $x$ can be rewritten as the combination
$$
x=\frac{j}\beta + \frac{z}\beta\,, \qquad
z= \sum_{i=1}^{k-1}\frac{x_{i+1}}{\beta^i}\,.
$$
By induction hypothesis, $z$ could be generated from $0$ and $1$ using
the given operations; i.e. $z\in{\it Cl}_{\mathcal M}\{0,1\}$. Since
\begin{equation}\label{eq:afin}
(x\Vdash_i y) + c = (x+c)\Vdash_i (y+c)\,,\quad\hbox{and }\quad
\frac{x}\beta \Vdash_i \frac{y}\beta = \frac1\beta (x\Vdash_i y)\,,
\qquad\hbox{for any } i=1,\dots,m\,,
\end{equation}
one has the following relation,
$$
{\it Cl}_{\mathcal M}\left\{\frac{j}\beta , \frac{j+1}\beta\right\}
= \frac{j}\beta + \frac1\beta {\it Cl}_{\mathcal M}\{0,1\}\,.
$$
Therefore
$$
x \,=\, \frac{j}\beta + \frac{z}\beta  \;\in\; 
\frac{j}\beta + \frac1\beta {\it Cl}_{\mathcal M}\{0,1\} \,=\, 
{\it Cl}_{\mathcal M}\left\{\frac{j}\beta , \frac{j+1}\beta\right\}
\subset {\it Cl}_{\mathcal M}\{0,1\}\,.
$$
The later inclusion is valid due to the first step of induction 
(see~\eqref{eq:1step}).

Assume that the coefficient $x_1$ of $x$ in its $\beta$-expansion 
is equal to $m$. Necessarily, from the properties of $\beta$-expansions,
$x_2$ is equal to 0. In this case we use 
$$
x \,=\, \frac{m}\beta + \frac{z}{\beta^2}  \;\in\; 
\frac{m}\beta + \frac1{\beta^2} {\it Cl}_{\mathcal M}\{0,1\} \,=\, 
{\it Cl}_{\mathcal M} \left\{ \frac{m}\beta , 1 \right\} \subset 
{\it Cl}_{\mathcal M}\{0,1\}\,.
$$

Thus we have shown that all points in $[0,1]\cap\Zb$ can be generated 
starting from 0, and 1, using the operations $\Vdash_i$, $i=1,\dots,m$.
The following relation can be used to reduce the number of necessary 
operations from $m$ to $[(m+1)/2]$. The relation is valid for any $k$, 
$k=1,\dots,m-1$, independently of $x,y$. 
$$
(x\Vdash_{k+1} y) \Vdash_1 (x\Vdash_k y) = 
\bigl(1-\frac{m-k}{\beta} \bigr) x + 
\frac{m-k}{\beta}  y = y \Vdash_{m-k} x\,.
$$

\item{(b)} 
Let now $\beta$ be the root of $x^2=mx-1$, $m\geq4$. In this case, 
if we want to use the $\beta$-expansions of numbers, we have to 
encounter a complication: There exist points in the ring $\Zb$ which 
do not have a finite $\beta$-expansion.

Recall that an $x\in\Zb$ has a finite $\beta$-expansion 
if and only if $N(x)=xx'\geq0$, (see~\eqref{eq:burdia}). For us, 
only $x\in(0,1)$ are of interest. The following statement is valid.
Let $x\in\Zb\cap(0,1)$. Then either $x$ or $(1-x)$ has a finite
$\beta$-expansion. Indeed, since $N(x)=xx'$ is an integer and $x<1$, 
necessarily $|x'|>1$. Therefore $x'>0$ if and only if $1-x'<0$.

Now we can proceed, similarly as for (a)
by induction on the length of the $\beta$-expansion in order to show that
any element of ${\rm Fin}(\beta)\cap[0,1]$ can be generated by corresponding
operations $\Vdash_i$, $i=1,\dots,m-1$.
For elements $x\in(0,1)\cap\bigl(\Zb\setminus{\rm Fin}(\beta)\bigr)$, 
we have $1-x\in(0,1)\cap{\rm Fin}(\beta)$. Therefore $1-x$ can be generated from
0 and 1 using the given operations.
In the corresponding combination we replace all 0 by 1 and vice versa, which
gives the desired combination for element $x$. 

The coefficients in a $\beta$-expansion take one of the values $0,\dots,m-1$.
For an $x\in{\rm Fin}(\beta)\cap(0,1)$, whose $x_1=j$, $j=0,\dots,m-2$, the
procedure is the same as in the case (a). Consider an
$$
x=\frac{m-1}\beta + \frac{x_2}{\beta^2} + \sum_{i=3}^k \frac{x_i}{\beta^i}
=1-\frac1\beta + \frac{x_2+1}{\beta^2} + \sum_{i=3}^k \frac{x_i}{\beta^i}
=\frac1\beta z + \left(1-\frac1\beta\right)1=z\Vdash_1 1 \,.
$$
Since the segment $x_1x_{2}\dots x_{k}$ of a $\beta$-expansion is 
strictly lexicographically smaller than the sequence $(m-1)(m-2)(m-2)\dots$, 
the string $x_2x_3\dots x_k$ is strictly smaller than $(m-2)(m-2)\dots$. 
Therefore  $(x_2+1)x_3\dots x_k$ is the $\beta$-expansion of $z$. 
By induction hypothesis $z\in{\it Cl}_{\mathcal M}\{0,1\}$, 
therefore $x=z\Vdash_1 1\in{\it Cl}_{\mathcal M}\{0,1\}$.

By that we have shown that any $x\in[0,1]\cap\Zb$ can be generated by 
given operations $\Vdash_i$, $i=1,\dots,m-1$. The number of operations
can be reduced from $m-1$ to $[(m-1)/2]$ by the following relation, which 
is valid for any $k$, $k=1,\dots,m-2$, for all $x,y$. 
$$
(x\Vdash_k y) \Vdash_1 (x\Vdash_{k+1} y) = 
\left(1-\frac{m-k-1}{\beta} \right) x + 
\frac{m-k-1}{\beta}  y = y \Vdash_{m-k-1} x\,.
$$
\end{trivlist}
\vskip-0.2in
\pfk

Due to~\eqref{eq:afin}, the Lemma~\ref{eins} can be generalized in 
the following way. If $c\in\Zb$, we have 
\begin{equation}\label{eq:skoky}
{\it Cl}_{\mathcal M}\{c,c+1\} = [c,c+1]\cap\Zb\,.
\end{equation}

\pf[Proof of Proposition~\ref{jajaj}]
In order to prove the statement of the proposition, we have to find 
$\Omega$, such that $\Lambda'=\Omega\cap\Zb$. We show that 
this property is satisfied by the convex hull of $\Lambda'$, so we 
put $\Omega=\langle \Lambda' \rangle$. We have to justify that any 
element of $\Zb\cap\langle\Lambda'\rangle$ can be generated using 
the given operations starting from points of $\Lambda'$. 

Consider $x\in\Zb\cap\langle\Lambda'\rangle$.
Since $0,1\in\Lambda$, then also $0,1\in\Lambda'$, and hence 
$\Lambda'\supset([0,1]\cap\Zb)$, as a consequence of Lemma~\ref{eins}.
Therefore, if $x\in[0,1]$, then $x$ belongs to $\Lambda'$.
Without loss of generality let $x>1$. 
Observe that the points of $\Lambda'$ cover densely the interval 
$\langle\Lambda'\rangle$. Therefore one can
find a finite sequence of intervals $[y_1-1,y_1],\dots,[y_k-1,y_k]$, 
such that the elements $y_i$ belong to $\Lambda'$, and it holds 
that $0<y_1-1<1$, $y_{i+1}-1<y_i$ and $y_k-1<x\leq y_k$.
According to~\eqref{eq:skoky}, one has 
$$
\Lambda'\supset{\it Cl}_{\mathcal M}\{y_i-1,y_i\}=[y_i-1,y_i]\cap\Zb\,.
$$
Therefore $x$ can be generated from points of $\Lambda'$ by the 
operations $\Vdash_i$, hence it is contained in $\Lambda'$.
Thus we have proven that $\Lambda'=\Omega\cap\Zb$.

It suffices now to justify that $\Omega$ is bounded. This is true, 
because otherwise $\Lambda=\Sigma_\beta(\Omega)$ would
not be uniformly discrete.
\pfk

\section{Main results, comments and open problems}\label{sec:results}

Due to Proposition~\ref{vylouc}, for any quadratic
Pisot number $\beta\neq\frac12(1+\sqrt5)$, $1+\sqrt2$,
$2+\sqrt3$, $s$-convexity of a uniformly discrete set
$\Lambda$ with an arbitrary fixed parameter $s$ does 
not mean that $\Lambda$ is a model set in $\Zb$. 
The example of an $s$-convex set which is not 
a model set in $\Zb$ is the closure ${\it Cl}_{s}\{0,1\}$ 
for arbitrary $s\in\Zb$, $s'\in(0,1)$. 

In the three exceptional cases $\beta=\frac12(1+\sqrt5)$, ($x^2=x+1$),
$\beta=1+\sqrt2$, ($x^2=2x+1$), $\beta=2+\sqrt3$, ($x^2=4x-1$),
the set ${\mathcal M}$ of parameters in Proposition~\ref{jajaj} has 
only one element $s'=1/\beta$, i.e.~$s=1/\beta'$. 
The $s$-convexity with this specific parameter $s$ implies being a model
set in the corresponding $\Zb$. However, Proposition~\ref{vylouc}
states that $s=1/\beta'$ and $1-s$ are the only parameters in $\Zb$ with
this property.

\begin{thm}\label{thm:1}
Let $\beta$ be a quadratic unitary Pisot number and $s\in\Zb$.
We say that a parameter $s$ is model-set-forcing if 
any uniformly discrete $s$-convex set $\Lambda\subset\Zb$, 
containing 0 and 1, is a model set in $\Zb$. 

The parameter $s$ is model-set-forcing if and only if $s$ or $1-s$
takes one of the three values $-\frac12(1+\sqrt5)$, $-1-\sqrt2$, 
or $2+\sqrt3$.
\end{thm}

In the majority of cases a single parameter $s\in\Zb$ does not suffice
to characterize model sets in $\Zb$ by $s$-convexity. The following 
theorem, which is a direct consequence of Proposition~\ref{jajaj}, states
that for such characterization one may use more, but a finite number
of parameters $s$. 

\begin{thm}\label{thm:2}
Let $\beta$ be quadratic unitary Pisot number, satisfying the equation 
$\beta^2=m\beta\pm1$. Let $\Lambda$ be a uniformly discrete subset 
of $\Zb$, containing 0,1. Set 
$$
{\mathcal N}=\left\{s \biggm | s'=\frac1\beta\,,\,\frac2\beta\,,\,\dots\,,\,
\frac1\beta\left[\frac{m\pm1}2\right]\right\}\,.
$$
Then $\Lambda$ is a model set in $\Zb$ if and only if $\Lambda$ is 
$s$-convex for any $s\in {\mathcal N}$.
\end{thm}

The only cases, when the set ${\mathcal N}$ has one 
element, are the three exceptional quadratic unitary Pisot numbers 
singled out by Theorem~\ref{thm:1}. For other irrationalities $\beta$, 
the two theorems state that for a characterization of model sets in 
$\Zb$ one has to consider at most $[(m\pm1)/2]$ but at least two 
parameters. It would be interesting to determine in which cases 
one may reduce the number of necessary parameters from the 
suggested $[(m\pm1)/2]$.

Note that if instead of a model set $\Sigma_\beta(\Omega)$ we consider
its multiple $\xi\Sigma_\beta(\Omega)$, for $\xi\in\Zb$, the resulting
set is generally not a model set in $\Zb$ according to the 
the definition~\eqref{eq:deqc}, but only its affine image.
In order to avoid such trivial examples of a uniformly discrete $s$-convex 
set $\Lambda$, which is not a model set, we require for our statements 
$0,1\in\Lambda$.
One may impose on $\Lambda$ a weaker assumption using the notion
of greatest common divisors in $\Zb$. However, 
the formulations and proofs would require complicated technicalities
without really obtaining something new.

Non-trivial examples of uniformly discrete $s$-convex sets, which are not
model sets in any $\Zb$, are the closures ${\it Cl}_s\{0,1\}$ for practically
all parameters $s\in\Zb$. Such closures are proper subsets of model sets.
However, it remains an open question, whether these subsets satisfy
in addition to the uniform discreteness also the other part of Delone property, 
namely the relative density. If ${\it Cl}_s\{0,1\}$ is Delone, it is an interesting
new aperiodic structure with abundant self-similarities. Investigation
of even basic properties of such sets would be desirable.

The content of the article concerns the relation of $s$-convexity and
1-dimensional model sets. Probably a more interesting question is, 
whether the notion of $s$-convexity can be used for a 
characterization of model sets in any dimensions. 
Such question was solved in~\cite{selfs} for model sets based on the 
golden mean irrationality. We have shown there that any $(-\tau)$-convex 
uniformly discrete set $\Lambda\subset\R^n$ is an affine image of a model 
set. We expect similar fact to be true for other quadratic unitary Pisot 
numbers as well.

\section*{Acknowledgements}

We are grateful to Drs.~Frank~W.~Lemire and Alfred~Weiss for useful
suggestions and comments.
We acknowledge partial support from the National Science and
Engineering Research Council of Canada, FCAR of Quebec, by
FR VS 767/1999 and by NATO Collaborative Research Grant CRG 974230.
Two of the authors (Z.M. and E.P.) are grateful for the hospitality of the 
Centre de recherches math\'ematiques, Universit\'e de Montr\'eal, 
and of the Joint Institute for Nuclear Research 
in Dubna. One of us (J.P.) is grateful for the hospitality of the
Mathematical Sciences Research Institute in Berkeley.



\end{document}